\documentclass[12pt]{article}
\usepackage{amssymb,amsmath,amsfonts,amsthm}
\usepackage{array}
\usepackage{bm}

\def\x#1{}                                                      
\def\eq#1{\begin{equation}#1\end{equation}}                     
\def\eqs*#1{\begin{eqnarray*}#1\end{eqnarray*}}                 
\def\eqss#1{\begin{eqnarray}#1\end{eqnarray}}                   
\newtheorem{defin}{Definition}{\bfseries}{\upshape}             
\def\Up#1{\vspace{-#1em}}                                       
\def\xy{\hspace{.07em}}                                         
\def\xz{\hspace{-.07em}}                                        
\def\a{\alpha}                                                  
\def\R{{\mathbb R}}                                             
\def\f{(\cdot,\cdot)}                                           
\def\minl{\mathop{\inf}\limits}                                 
\def\maxl{\mathop{\sup}\limits}                                 
\def\RL{\,\Leftrightarrow\,}                                    
\def\s{\sigma}                                                  
\def\SS{\Sigma}                                                 
\def\gee{\xy\ge\xy}                                             

\author{     Michel Deza\footnotemark[1]
   \and Pavel Chebotarev\footnotemark[4]}

\title{\vspace{-2.5em}Protometrics\\
      {\vspace{-0.4em}\small\emph{In honor of Ivo Rosenberg}}\Up{.5}}
\date{}

\def\thefootnote{\fnsymbol{footnote}}

\begin{document}

\footnotetext[1]{Ecole Normale Superieure, Paris,                                          e-mail:\;{\footnotesize\tt Michel.Deza@ens.fr.}}
\footnotetext[4]{Institute of Control Sciences of the Russian Academy of Sciences, Moscow, e-mail:\:{\footnotesize\tt pavel4e@gmail.com.}}\x{2016}
\maketitle
\vspace{-2em}
\begin{abstract}
\noindent
We introduce the concept of protometric 
and present some properties of protometrics.
\end{abstract}
\def\thefootnote{\arabic{footnote}}

This note is a tribute to Ivo, a friend and co-author of the first author during last 40 years.
It is written in the taste and style of Ivo, on the border of Logic, Distance Spaces and Combinatorics.
\begin{defin}
\label{d_3}
{\rm
For a set $X,$ we say that a function $d\xz\!:\!X\!\times\xz X\!\to\!\R$ satisfies the {\em triangle inequality of type\/}:
$$
\begin{tabular}{ll>$l<$}
o (outgoing)                   &iff &d(x,y)+d(x,z)\gee d(y,z);\\
i (incoming)\footnotemark  [1] &iff &d(y,x)+d(z,x)\gee d(y,z);\\
t (transitive)\footnotemark[2] &iff &d(y,x)+d(x,z)\gee d(y,z);\\
c (cyclic)                     &iff &d(z,x)+d(x,y)\gee d(y,z)
\end{tabular}
$$
for all ${\,x,y,z\in X}.$ These will also be termed *-triangle inequalities, where * are type letters.
}
\end{defin}
\footnotetext[1]{It is also called the \emph{strong   triangle inequality.}}
\footnotetext[2]{It is also called the \emph{oriented triangle inequality\/} or simply the \emph{triangle inequality.}}
\setcounter{footnote}{2}

\noindent
{\bf Simple facts about the triangle inequalities}
\begin{enumerate}
\Up{.4}\item
{\em If $d\f$ is symmetric$,$ i.e.$,$ $d(x,y)=d(y,x)$ for all $x,y\in X,$ then the four versions of the triangle inequality are equivalent.}

\Up{.4}\item
{\em Let $d'(x,y)=d(y,x)$ for all $x,y\in X.$
Then $d'\f$ satisfies the i-triangle inequality iff $d\f$ satisfies the o-triangle inequality $(\!$and vice versa$)$.
The t-triangle inequality holds or fails for $d\f$ and $\,d'\f$ simultaneously\/$,$ and so does the c-triangle inequality.}

\Up{.4}\item
{\em If $\,d\f$ satisfies the triangle inequality of type\/$:$

\Up{1.0}$$
\begin{tabular}{ll>$l<$}
o\/$,$  &then &d(x,x)\in\Bigl[\,\maxl_{y\in X}\xy(d(y,x)-d(x,y)),\,2\minl_{y\in X}    d(y,x)\Bigr];\\
i\/$,$  &then &d(x,x)\in\Bigl[\,\maxl_{y\in X}\xy(d(x,y)-d(y,x)),\,2\minl_{y\in X}    d(x,y)\Bigr];\\
t\/$,$  &then &d(x,x)\in\Bigl[\xy0,\,                               \minl_{y\in X}\xy(d(x,y)+d(y,x))\Bigr];\\
c\/$,$  &then &d(x,x)\in\Bigl[\,\maxl_{y\in X}|d(x,y)-d(y,x)|,\,    \minl_{y\in X}\xy(d(x,y)+d(y,x))\Bigr]
\end{tabular}
$$
for all ${\,x\!\in\! X}.$ Hence\/\footnote{Moreover, for types t or c, $\,d(x,y)\ge0\xy$ holds if the symmetry of $d\f$ is additionally assumed.}
$\xy d(x,x)\ge0\xy$ for types~t or~c and $\,d(x,y)\ge0\xy$ for types~o or~i.
}

\Up{.4}\item
{\em If $\,[\xy d(x,y)\!=\!0\RL x\!=\!y\,]\xy$ \emph{(identity of indiscernibles)} and $d\f$ satisfies the triangle inequality of type~o\/$,$ i\/$,$ or c\/$,$ then $d\f$ is symmetric and nonnegative\/$,$ thus\/$,$ $d\f$ is a metric.}
\end{enumerate}

\begin{defin}
\label{d_2}
{\rm
For a set $X,$ a function $p\xz\!:\!X\!\times\xz X\!\to\!\R$ satisfies the {\em pre-quadrangle inequality of type\/}:
$$
\begin{tabular}{ll>$l<$}
o (outgoing)\footnotemark[4]   &iff &p(x,y)+p(x,z)\gee p(y,z)+p(x,x);\\
i (incoming)                   &iff &p(y,x)+p(z,x)\gee p(y,z)+p(x,x);\\
t (transitive)\footnotemark[5] &iff &p(y,x)+p(x,z)\gee p(y,z)+p(x,x);\\
c (cyclic)                     &iff &p(z,x)+p(x,y)\gee p(y,z)+p(x,x)
\end{tabular}
$$
for all ${x,y,z\!\in\! X}.$
Such a function $p\f$ is called a {\em protometric\/} of the corresponding type.
If~the inequality is strict whenever $z=y$ and $y\ne x$ \cite{CheSha98S}, then $p\f$ is a \emph{strict protometric.}}
\end{defin}
\footnotetext[4]{In this case, $-p\f$\x{2016} satisfies the \emph{triangle inequality for proximities\/} appeared in \cite{CheSha98S} and some earlier papers by the same authors; cf.\:\cite{CatralNeumannXu05}.}\x{}
\footnotetext[5]{This inequality was considered by Matthews~\cite{Matt92}. It is also called the \emph{sharp triangle inequality} and the \emph{modified triangle inequation.}\x{}
\setcounter{footnote}{4}}

\noindent
{\bf Simple facts about protometrics}
\begin{enumerate}
\Up{.4}\item\label{i_1}
{\em The pre-quadrangle inequality of each type
\emph{strengthens} the triangle inequality of the same type if $\,[\,p(x,x)\ge0,$ but $\,p(x,x)\not\equiv0\xy\xy]$ and
\emph{reduces to it} whenever 
$\,p(x,x)\equiv0.$
Any metric is a protometric of each type.}

\Up{.4}\item\label{i_2}
{\em If $p\f$ is a protometric of type~o\/$,$ i\/$,$ or c\/$,$ then $p(x,y)\ge\frac12(p(x,x)+p(y,y))$ and $p\f$ is symmetric.
     If $p\f$ is a protometric of type~t\/$,$ then ${\,p(x,y)+p(y,x)\ge p(x,x)+p(y,y).}$
Thus$,$ there are only two types of protometrics\/$:$ general \emph{protometrics (\xz\xz}of type~t\/$)$ and 
\emph{symmetric protometrics (\!}of type o-i-t-c$)\xz$ also called\/$,$ in the case of $p(x,x)\ge0,$ \emph{weak partial pseudo-metrics~\cite{Heckmann99}.}}\x{}

\Up{.4}\item\label{i_3}
{\em 
For any $f\xz\!:\!X\!\to\!\R,$ $\xy p(x,y)\!=\!f(x)$\x{} and $\xy p(x,y)\!=\!f(y)$ are protometrics.
If $p$ is a proto\-metric\/$,$ then so is $p'(x,y)\!\stackrel{\rm def}=\!p(y,x).$
If $p$ and $q$ are protometrics on $X,$ then so is $p+q.$
Hence $p(x,y)+p(y,x)$ is a symmetric protometric whenever $p$ is a protometric.}

\Up{.4}\item\label{i_4}
{\em Let $\,p'(x,y)=\a\xy p(x,y)+f(x)+f(y)\,$ for all $\,x,y\in X,\,$ where $\a>0$ and $f\xz\!:\!X\!\to\!\R.$ Then $p$ and $\xy p'$ are or are not protometrics of the same type simultaneously.}

\Up{.4}\item\label{i_5}
{\em If in Fact\/~$\ref{i_4}$ it holds that $\,f(x)=-\frac\a2\xy p(x,x),$ then $p$\x{} is a protometric of any type iff $\,p'$\x{} satisfies the triangle inequality of the same type.} 

\Up{.4}\item\label{i_6}
{\em It follows from Facts~$\ref{i_2},$ $\ref{i_3},$ and~$\ref{i_5}$ that for any protometric $p,$ the symmetric function $\,d(x,y)=\a(p(x,y)+p(y,x)-p(x,x)-p(y,y)),$ $\a>0,$ satisfies the triangle inequality$,$ is non-negative$,$ and $\,d(x,x)\equiv0$.
If$,$ additionally\/$,$ $p$ is a strict protometric\/$,$ then ${\,d(x,y)>0}$ whenever $\xy y\ne x$ and thus\/$,$ $d$ is a metric.}
\end{enumerate}

\smallskip\noindent
{\bf Protometrics and similarity measures}

\smallskip
The concept of protometric as applied to dissimilarity measures $d,$\x{} is somewhat exotic, since $d(x,x)\!=\!0$ is characteristic of such measures.
However, the fulfillment of the pre-quadrangle inequality is quite typical of the function $\,-s\f,$ where $s\f$ is a similarity measure.

One example is the \emph{Gromov product similarity} (or \emph{covariance}) $(x.y)_{x_0}\!=\frac12(d(x,x_0)+d(y,x_0)-d(x,y)),$ where $d\f$ is a metric and $x_0\in X$ is any fixed \emph{base point}. It follows from the above Facts\;\ref{i_1} and\;\ref{i_4} that $\,-(x.y)_{x_0}$ is a non-positive symmetric protometric (cf.\ the {\em Farris transform metric\/} $\,C-(x.y)_{x_0},\,$ where $C$ is a large enough positive constant \cite[Chapter\:4]{DezaDeza09}).

Moreover, in \cite{CheSha98} a number of proximity measures $s\f$ for graph vertices were presented such that $\,-s\f$ are non-positive protometrics.

A function $s\xz\!:\!X\!\times\xz X\!\to\!\R$ satisfies the {\em transition inequality\/} if $s(y,x)\xy s(x,z)\le s(y,z)\xy s(x,x).$ In  \cite{Che11AAM}, it was shown that a number of positive proximity measures $s\f$ for the vertices of a strongly connected weighted digraph satisfy the transition inequality. Consequently, the corresponding functions $\,-\ln s\f$ are protometrics.

By Fact\;\ref{i_6}, for any strict protometric $p,$ $\,d(x,y)=p(x,y)+p(y,x)-p(x,x)-p(y,y)$ is a metric. In \cite{CheSha98S} it was shown that for the classes of $\SS_m$-proximities, this transformation is invertible. The \emph{$\SS_m$-proximities\/} are strict protometrics $\s$ such that for every $x\in X,$ ${\s(x,\cdot)=m,}$ where $\s(x,\cdot)$ is the value of an averaging linear functional applied to $\s(x,y)$ as a function of~$y.$ Alternatively, this result can be derived in terms of difference protometrics (see below).

\bigskip\noindent
{\bf $\bm0$-protometrics and difference protometrics}

\smallskip
$0$-protometrics are the antipodes of strict protometrics (Definition\;\ref{d_2}).
A~{\em $0$-protometric} is a protometric $p$ for which the pre-quadrangle inequality with $z\!=\!y$ always holds as equality: ${p(x,y)+p(y,x)-p(x,x)-p(y,y)\equiv0.}$
The $0$-protometrics form the largest linear space in the flat convex cone of protometrics on~$X.$ For a finite $X,$ a basis of this space is given by all but one $0$-protometrics $q'_u(x,y)\!=\!1_{x=u}$ and $q''_u(x,y)\!=\!1_{y=u}$ (since $\sum_uq'_u\equiv\sum_uq''_u\equiv1$). A~basis of the space of symmetric $0$-protometrics on $X$ is given by all $q'_u+q''_u$.

\smallskip
A \emph{difference protometric\/} (called \emph{strong protometric\/} in~\cite{DezaDeza09}) is a protometric $d$ such that $\,d(x,x)\equiv0.$
It satisfies $d(x,y)+d(y,x)\ge0.$ If, moreover, $d(x,y)\ge0$ for all $x,y\in X,$ then $d$\x{} is a \emph{quasi-semi-metric}~\cite[Section~1.1]{DezaDeza09}.
In this case, the binary relation $\xy x\preceq y\RL d(x,y)=0\xy$ is a preorder on~$X$. Note that it is the \emph{specialization preorder\/} in the induced topology.

Given a difference protometric $d$ on $X$ and any function $f\xz\!:\!X\!\to\!\R,$
\eq{
\label{e_de}
p(x,y)=\frac12(d(x,y)+f(x)+f(y))
}
defines a protometric by virtue of Fact~\ref{i_4}. It follows from \eqref{e_de} that for any $x,y\in X,$
\eqss{
\label{e_dif}
f(x)&=& p(x,x);\\
\label{e_di}
d(x,y)&=&2p(x,y)-p(x,x)-p(y,y).
}
Thus, mappings \eqref{e_de} and \eqref{e_dif}--\eqref{e_di} establish a bijection between all pairs $(d,f),$ where $d$ is a difference protometric on $X,\,$ $f$ being a function $X\!\to\!\R,$ and some protometrics $p$ on~$X.$ Since each protometric $p$ generates a difference protometric $d$ by means of \eqref{e_di}, this bijection involves all protometrics $p$ on~$X.$ The same mappings establish a bijection between the pairs $(d,f),$ where $d$ is a semi-metric, and the symmetric protometrics on~$X.$ For example, a pair $(d,f)$ with $d$ being a semi-metric and $f(x)=-d(x,x_0)$ corresponds to the protometric $-(x.y)_{x_{0}}$.
Since $f(x)\xz+\xz f(y)$ is a 0-protometric, \eqref{e_de} in the symmetric case with a finite $X$ enables the representation
$$
p=\frac12\Bigl(d+\sum_{u\in X}p(u,u)(q'_u+q''_u)\!\Bigr).
$$

The \emph{difference $0$-protometrics}\x{2016} (the elements of the intersection of difference protometrics and $0$-protometrics) are exactly \emph{potential differences\/}, i.e., functions of the form $d(x,y)=h(x)-h(y),$ where $h\xz\!:\!X\!\to\!\R.$

\smallskip
Note that the protometrics of this paper are not related to the \emph{proto-metrizable}, i.e., paracompact and having an orthobase, topological spaces.

\end{document}